\documentclass[11pt]{amsart}
\usepackage{amsmath}
\usepackage{amssymb}
\usepackage{amscd}
\usepackage{url}
\usepackage[latin2]{inputenc}
\usepackage[T1]{fontenc}
\usepackage{pstricks}
\usepackage{pst-node}
\usepackage{pst-coil}
\usepackage{graphicx}

\newcommand{\intfrac}[2]{\genfrac{\lfloor}{\rfloor}{}{1}{#1}{#2}}

\newcommand{\intpart}[1]{\left\lfloor#1\right\rfloor}
\newcommand{\sawtooth}[1]{\left\langle #1\right\rangle}

\def\R{{\mathbb R}}
\def\spq{{\sum_{s=0}^{pq-1}}}
\def\soi{{\sum_{k=0}^{p-1}}}
\def\si{{\sum_{k=1}^{p-1}}}
\def\soj{{\sum_{l=0}^{q-1}}}
\def\sj{{\sum_{l=1}^{q-1}}}
\def\sord{\sigma_{ord}}

\theoremstyle{definition}
\newtheorem{example}{Example}[section]

\theoremstyle{plain}

\newtheorem{lemma}[example]{Lemma}
\newtheorem{proposition}[example]{Proposition}

\newtheorem{corollary}[example]{Corollary}

\theoremstyle{remark}
\newtheorem*{acknowledgements}{Acknowledgements}
\newtheorem{remark}[example]{Remark}
\numberwithin{equation}{section}

\newcommand{\jed}{\mathbf{1}}
\makeatletter
\@namedef{subjclassname@2010}{\textup{2010} Mathematics Subject Classification}
\makeatother
\title[On the signatures of torus knots]{On the signatures of torus knots}
\author{Maciej Borodzik}
\address{Institute of Mathematics, University of Warsaw, ul. Banacha 2,
02-097 Warsaw, Poland}
\email{mcboro@mimuw.edu.pl}
\thanks{The first author is supported by Polish MNiSzW Grant No N N201 397937 and also by the Foundation for Polish Science FNP.
The second author is supported by Polish MNiSzW Grant N N201 397437.}
\author{Krzysztof Oleszkiewicz}
\address{Institute of Mathematics, University of Warsaw, ul. Banacha 2,
02-097 Warsaw, Poland;
Institute of Mathematics, Polish Academy of Sciences, ul. \'Sniadeckich 8, 00-956 Warsaw, Poland}
\email{koles@mimuw.edu.pl}

\date{24 February 2010}

\begin{document}
\begin{abstract}
We study properties of the signature function of the torus knot
$T_{p,q}$. First we provide a very elementary proof of the formula
for the integral of the signatures over the circle. We obtain also
a closed formula for the Tristram--Levine signature of a torus
knot in terms of Dedekind sums.
\end{abstract}
\maketitle

\section{Preliminaries}

Let $K$ be a knot in $S^3$ with a Seifert matrix $S$. Let also $z\in S^1$, $z\neq 1$ be a complex number.
The \emph{Tristram--Levine} signature $\sigma(z)$ is the signature of the hermitian form
\[(1-z)S+(1-\bar{z})S^T.\]
This is obviously an integer-valued piecewise constant function.
It does not depend on a particular choice of Seifert matrix. If we
substitute $z=-1$ we get an invariant $\sord$, which is called the
\emph{(ordinary) signature}. We define also the integral $I_K$
\[I_K=\int_0^1\sigma(e^{2\pi ix})\,dx.\]

Signatures are very strong knot cobordism invariants, which can be
used to bound the four-genus and the unknotting number of $K$. The
integral $I_K$ of the signature function is one of the so called
$\rho$ invariants of knots (see \cite{COT1,COT2}) and is 
of independent interest.

For a torus knot $T_{p,q}$, where $\gcd(p,q)=1$, the signature
function can be expressed in the following nice way (see \cite{Li}
or \cite[Chapter XII]{Kau})
\begin{proposition}\label{p:litherland}
Let
\begin{equation}\label{eq:sigmadef}
\Sigma=\left\{\dfrac{k}{p}+\dfrac{l}{q}\colon 1\le k\le p-1,\,\,1\le l\le q-1\right\}.
\end{equation}
Then for any $x\in(0,1)\setminus\Sigma$ we have
\begin{equation}\label{eq:sig}
\sigma(e^{2\pi ix})=|\Sigma\setminus(x,x+1)|-|\Sigma\cap(x,x+1)|,
\end{equation}
where $|\cdot|$ denotes the cardinality of a set.
In particular $\sord=|\Sigma\setminus(1/2,3/2)|-|\Sigma\cap(1/2,3/2)|$.
\end{proposition}

The explicit formulae for $\sord$ and $I_K$ of torus knots have
been known in the literature for quite a long time. In fact,
$\sord$ by a result of Viro (see \eqref{eq:tau2}) is equal 
to $\tau_2$, which was computed in \cite{HZ} for $p$ and $q$ odd, 
and (denoted as $\sigma(f+z^2)$) in \cite{Nem} in general case. 
On the other hand, Kirby and Melvin
\cite[Remark 3.9]{KM} and \cite[Example 4.3]{Nem} provided a
formula for $I_K$. Nevertheless all the above-mentioned results
are related more to singularity theory and low-dimensional
topology, than to knot theory itself.

After the discovery of $\rho$ invariants, the interest of computing
$I_K$ for various families of knots grew significantly. Two
independent new proofs of the formula for $I_K$ of torus knots
\cite{Bo,Co} appeared in 2009. In particular \cite{Bo} provided a
bridge between the $I_K$ and cuspidal singularities of plane
curves.

In this paper we present an elementary proof of the formula for
$I_K$ (Proposition~\ref{prop:int}). We also cite a formula of
N\'emethi and draw some consequences from it. In
Section~\ref{S:val} we use a theorem of Rosen to obtain the
explicit value of the signature $\sigma(z)$ of a torus knot not
only for $z=-1$, but also for any $z\in S^1\setminus\{1\}$ (Proposition~\ref{prop:siggen}). 
This result seems to be new. In
Section~\ref{S:proofnogo} we show that the formula for
$\sord(T_{p,q})$ cannot be written as a rational function of $p$
and $q$.

\section{Formula for the integral}\label{S:int}
\begin{proposition}\label{prop:int}
For a torus knot $T_{p,q}$ we have
\begin{equation}\label{eq:I}
I=-\dfrac13\left(p-\dfrac1p\right)\left(q-\dfrac1q\right).
\end{equation}
\end{proposition}
This proposition was first proved in \cite[Remark 3.9]{KM}. Refer to \cite{Nem,Bo,Co} for other proofs.

\begin{proof}
Let $f(x)=-\sigma(e^{2\pi i x})$ and $J=\int_0^1f(x)\,dx=-I$. Then
\[f(x)=\sum_{y\in\Sigma}\jed_{(y,y+1)}(x)-\sum_{y \in \Sigma} \jed_{\R\setminus (y,y+1)}(x).\]
(Here, for a set $A\subset\mathbb{R}$, $\jed_A$ denotes the function which is equal to $1$ on $A$ and $0$ away from $A$.)
Hence
\[J=\sum_{y \in \Sigma} \int_{0}^{1} \left( \jed_{(y-1,y)}(x)-\jed_{\R \setminus (y-1,y)}(x)\right)\,dx=
\sum_{y \in \Sigma} (1-2|y-1|).\]
It follows that
\[J=\si \sj \left(1-2\left| \dfrac{k}{p}+\dfrac{l}{q}-1 \right|\right).\]
As for any $u,v\in\R$ we have $1-2|u+v-1|=2\min(1-u,v)+2\min(u,1-v)-1$,
\begin{multline*}
J=2\si \sj \min\left( \dfrac{p-k}{p}, \dfrac{l}{q} \right) +2\si \sj \min\left( \dfrac{k}{p}, \dfrac{q-l}{q} \right) -(p-1)(q-1)=\\
=4\si \sj \min\Big( \dfrac{k}{p}, \dfrac{l}{q} \Big)-(p-1)(q-1)=\dfrac{4}{pq}\si \sj \min(qk,pl)-(p-1)(q-1).
\end{multline*}
Now, obviously,
\begin{multline*}
\si \sj \min(qk,pl) = \\
=\sum_{s=0}^{\infty} \left|\left\{ \{ 1,\ldots, p-1\} \times \{ 1, \ldots, q-1\}\colon qk>s \text{ and } pl>s \right\}\right|=\\
=\spq (p-1-\intpart{s/q})(q-1-\intpart{s/p}).
\end{multline*}

We can multiply the expression in parentheses. Then, as $\spq \intpart{s/p}=p\soj l=\dfrac{1}{2}pq(q-1)$ we get

\begin{multline*}
\spq (p-1-\intpart{s/q})(q-1-\intpart{s/p})=pq(p-1)(q-1)-\dfrac{1}{2}pq(p-1)(q-1)-\\ \dfrac{1}{2}pq(p-1)(q-1)+
\spq \intpart{s/p}\intpart{s/q}=\spq \intpart{s/p}\intpart{s/q}.
\end{multline*}

It remains to compute $\spq \intpart{s/p}\intpart{s/q}$. To this
end let us denote by $R_{p}(s)$ the remainder of $s$ modulo $p$.
We then have

\begin{align*}
\spq \intpart{s/p}\intpart{s/q}= \spq \left( \dfrac{s-R_{p}(s)}{p} \cdot \dfrac{s-R_{q}}{q} \right)=\\
\dfrac{1}{pq}\left( \spq s^{2}- \spq sR_{p}(s) -\spq sR_{q}(s)+ \spq R_{p}(s)R_{q}(s)\right) =\\
\dfrac{1}{3}p^{2}q^{2}+\dfrac{1}{4}pq-\dfrac{1}{4}p^{2}q-\dfrac{1}{4}pq^{2}-\dfrac{1}{12}p^{2}-\dfrac{1}{12}q^{2}+\dfrac{1}{12},
\end{align*}
where we used the fact that $\spq R_{p}(s)R_{q}(s)=\soi \soj kl$ by the Chinese remainder theorem.

Putting all the pieces together we obtain the desired formula.
\end{proof}

Let us now present another proof, due to N\'emethi \cite{Nem}, see also \cite{Br,HZ}. Before
we do this, let us recall some facts from topology.

Assume that the knot $K$ is drawn on $S^3=\partial B^4$ and consider a Seifert surface $F$ of $K$. Let us push it
slightly into $B^4$ and for an integer $m$ let $N_m$ be the $m$ fold cyclic cover of $B^4$ branched along $F$.
Then the quantity $\tau_m=\sigma(N_m)$ (here $\sigma$ is a signature of a four-manifold)
is independent of the choices made. We have the
formula essentially due to Viro (see \cite[Section 2]{GLM} or \cite{Vi}).

\begin{equation}\label{eq:taum}\tau_m=\sum_{k=1}^{m-1}\sigma_K(\xi^k),
\end{equation}

where $\xi$ is a primitive root of unity of order $m$. In particular, since $\sigma$ is a Riemann integrable function, we have

\begin{equation}\label{eq:limit}
I=\int_0^1\sigma(e^{2\pi i x}),dx=\lim_{m\to\infty}\dfrac1m\tau_m.
\end{equation}

On the other hand

\begin{equation}\label{eq:tau2}
\tau_2(K)=\sord(K).
\end{equation}

If $K$ is a torus knot $T_{p,q}$ and $m,p,q$ are pairwise coprime,
then the $m$-fold cover of $S^3$ branched along $K$ is diffeomorphic to the Brieskorn
homology sphere $B(p,q,m)$ (see \cite{Br}, \cite[Section 5]{GLM}). Then $\tau_m$ turns
out \cite[Section 10.2 and 11]{HZ}
to be the signature of the manifold $X_{p,q,m}$ defined as the intersection of $z_1^p+z_2^q+z_3^m=\varepsilon$ with
$B(0,1)\subset \mathbb{C}^3$. In this context $\tau_m$ was computed by \cite[Formula 11 on page 122]{HZ} and
by \cite[Example 4.3]{Nem}. Especially the last formula is worth citing (N\'emethi uses $m(S(f))$ to denote the limit
\eqref{eq:limit}).

\begin{equation}\label{eq:nem}
I=-4(s(p,q)+s(q,p)+s(1,pq)).
\end{equation}

Here $s(a,b)$ is the Dedekind sum (see Section~\ref{S:ros}).
As by elementary computations $s(1,pq)=\dfrac{(pq-1)(pq-2)}{12pq}$, we get that
\[s(p,q)+s(q,p)=-\dfrac{I}{4}-\dfrac{(pq-1)(pq-2)}{12pq}.\]
Now we can look at the above equation as defining $I$ in terms of
$s(p,q)+s(q,p)$, but if we know $I$, we know $s(p,q)+s(q,p)$. In
other words we get the following observation.

\begin{corollary}
Any proof of Proposition~\ref{prop:int} provides a proof of the Dedekind reciprocity law.
\end{corollary}

\section{Lattice points in the triangle}\label{S:ros}

Let us recall basic definitions. For a real number $x$, $\intpart{x}$ denotes the integer part and $\{x\}=x-\intpart{x}$
the fractional part. The \emph{sawtooth} function is defined as

\[\sawtooth{x}=\begin{cases} \{x\}-\dfrac12& x\not\in\mathbb{Z}\\0& x\in\mathbb{Z}.\end{cases}\]

Sometimes $\sawtooth{x}$ is denoted $((x))$. We prefer this notation because it does not lead to confusion
with ordinary parenthesis.
We can now define the functions (below $p$, $q$ and $m$ are integers and $x,y$ are real numbers):

\begin{align*}
s(p,q)&=\sum_{j=0}^{p-1}\sawtooth{\dfrac{j}{q}}\sawtooth{\dfrac{pj}{q}}\\
s(p,q;x,y)&=\sum_{j=0}^{p-1}\sawtooth{\dfrac{j+y}{q}}\sawtooth{p\dfrac{j+y}{q}+x}.
\end{align*}

These functions satisfy the following reciprocity laws (see \cite{RG,HZ}). 
If $m$, $p$ and $q$ are pairwise coprime, then

\begin{align}
s(p,q)+s(q,p)&=\dfrac{1}{12}\left(\dfrac{p}{q}+\dfrac{q}{p}+\dfrac{1}{pq}\right)-\dfrac14\label{eq:recips}\\
s(p,q,x,y)+s(q,p,y,x)&=-\dfrac14d(x)d(y)+\sawtooth{x}\sawtooth{y}+\notag\\
+&\dfrac12\left(\dfrac{q}{p}\Psi_2(y)+\dfrac{1}{pq}\Psi_2(py+qx)+\dfrac{p}{q}\Psi_2(x)\right)\label{eq:recipss}
\end{align}

 Here

\[d(x)=\begin{cases} 1& \text{if $x\in\mathbb{Z}$}\\0&\text{otherwise}\end{cases}\]

and

\[\Psi_2(x)=B_2(\{x\})=\{x\}^2-\{x\}+\dfrac16\]

is the second Bernoulli polynomial. Now for a fixed $C\in[0,1)$ and $p$, $q$ coprime, let

\[A(p,q;C)=\{(k,l)\in\mathbb{Z}^2_{\ge 0}\colon 0\le \dfrac{k}{p}+\dfrac{l}{q}<1-C\}\]

and

\[N(p,q;C)=|A(p,q;C)|.\]

We have the following result due to Rosen \cite[Theorem~3.4]{Ro}.

\begin{proposition}\label{P:rosen}
In this case

\begin{equation}\label{eq:rosen}
\begin{split}
N(p,q;C)&=\dfrac{(1-C)^2}{2}pq+\dfrac{(1-C)}{2}(p+q)+\dfrac{q}{12p}+\dfrac{p}{12q}+K-\\
&-s(p,q;Cp,0)-s(q,p;Cq,0)+\sawtooth{Cp}+\sawtooth{Cq}+\\
&+(1-C)\sawtooth{Cpq}-(\dfrac78\delta_0+\dfrac38\delta_1-\dfrac18\delta_2)+\dfrac14,
\end{split}
\end{equation}

where

\[K=\begin{cases}\dfrac{1}{12pq}-\dfrac18&\text{if $Cpq\in\mathbb{Z}$}\\ \dfrac{1}{2pq}\Psi_2(Cpq)&\text{ otherwise}
\end{cases}\]

And for $r=0,1,2$, $\delta_r$ is the number of non-negative integers $k,l$ such that $\dfrac{k}{p}+\dfrac{l}{q}+C=r$.
\end{proposition}

This proposition admits an important corollary \cite[Corollary 3.5]{Ro}.

\begin{corollary}\label{c:rosen2}
If $p$ and $q$ are odd and coprime, then

\[N(p,q;\dfrac12)=\dfrac{pq}{8}+\dfrac{p+q}{4}+\dfrac{q}{6p}+\dfrac{p}{6q}+\dfrac{1}{24pq}-s(2p,q)-s(2q,p).\]

If $p$ and $q$ are coprime and $q$ is even, then

\begin{equation}\label{eq:rosen1/2}
N(p,q;\dfrac12)=\dfrac{pq}{8}+\dfrac{p+q}{4}-s(2p,q)+2s(p,q).
\end{equation}

\end{corollary}

We shall use these results to compute the signature of the torus knots. We need a following trivial lemma

\begin{lemma}\label{l:trivial}
The number of points $(k,l)\in A(p,q;C)$ such that $kl=0$ is equal to
\[Z(p,q;C)=\intpart{(1-C)p}+\intpart{(1-C)q}+1-d((1-C)p)-d((1-C)q),\]
where $d(x)$ again is $1$, if $x\in\mathbb{Z}$, and $0$ otherwise.
\end{lemma}

If $Cp$ and $Cq$ are not integers,
\[Z(p,q;C)=(1-C)(p+q)-\sawtooth{(1-C)p}-\sawtooth{(1-C)q}.\]

\section{Explicit formulae for the signatures}\label{S:val}
We begin with computing the value of the ordinary signature. As it
was already mentioned, $\sord=\tau_2$ (see \eqref{eq:tau2}) so the
first result below is in general known \cite{HZ,Nem}, but not
necessarily in the context of knot theory.

\begin{proposition}\label{th:one}
If $p$ and $q$ are both odd and coprime, then the ordinary signature of the torus knot $T_{p,q}$ satisfies

\[
\sord(T_{p,q})=-\dfrac{pq}{2}+\dfrac{2p}{3q}+\dfrac{2q}{3p}+\dfrac{1}{6pq}-4(s(2p,q)+s(2q,p))-1,
\]

where $s(x,y)$ is the Dedekind sum \emph{(}see Section~\ref{S:ros} or \cite{RG}\emph{)} (compare with
\cite[Formula 11 on page 122]{HZ}). If $p$ is
odd and $q>2$ is even, then
\[
\sord(T_{p,q})=-\dfrac{pq}{2}+1+4s(2p,q)-8s(p,q).
\]

\end{proposition}
\begin{proof}
Let us consider the torus knot $T_{p,q}$ and let $\Sigma$ be as in \eqref{eq:sigmadef}.
We can write  $\sord$ as
\begin{equation}\label{eq:S1}
\sord=4|\Sigma \cap (0,\dfrac12)|-|\Sigma|.
\end{equation}
Since $|\Sigma|=(p-1)(q-1)$, we need to find a closed formula for

\begin{equation}\label{eq:S}
S(p,q)=|\Sigma| \cap (0,\dfrac12)=\left|\left\{\dfrac{k}{p}+\dfrac{l}{q}<\dfrac12,\,\,1\le k\le p-1,\,\,1\le l\le q-1\right\}\right|.
\end{equation}

From the definition we get immediately that
\[S(p,q)=N(p,q;\dfrac12)-Z(p,q;\dfrac12).\]
Now $Z(p,q;\dfrac12)=\dfrac12(p+q)$ if $p$ and $q$ are both odd and $\dfrac12(p+q-1)$ if $q$ is even and $q>2$.
Hence, for $p$ and $q$ odd we have
\[S(p,q)=\dfrac{pq}{8}-\dfrac{p+q}{4}-s(2p,q)+2s(p,q),\]
while for $q$ even we have by \eqref{eq:rosen1/2}
\[S(p,q)=\dfrac{pq}{8}-\dfrac{p+q}{4}+\dfrac12-s(2p,q)+2s(p,q).\]
and using \eqref{eq:S1} we complete the proof.
\end{proof}

\smallskip
To express explicitly the values of Tristram--Levine signatures at
other points let us assume that $Cpq$ is not an integer. Define

\begin{align*}
M(p,q;C)&=N(p,q;C)-Z(p,q;C)=\dfrac{(1-C)^2}{2}pq-\dfrac{(1-C)}{2}(p+q)\\
&+\dfrac{q}{12p}+\dfrac{p}{12q}-s(p,q;Cp,0)-s(q,p;Cq,0)+\dfrac14-\\
&-\dfrac12(\sawtooth{Cp}+\sawtooth{Cq})+(1-C)\sawtooth{Cpq}+\dfrac{1}{2pq}\Psi_2(Cpq).
\end{align*}

Now it is a trivial consequence of Proposition~\ref{p:litherland} that if $C\in[0,1)$ and $e^{2\pi iC}=z$, then
\[\sigma(z)=-(p-1)(q-1)+2M(p,q;C)+2M(p,q;1-C).\]
Now, since for any integer $k$ and real $x$ we have $\sawtooth{(1-x)k}+\sawtooth{xk}=0$, the formula for
$M(p,q;C)+M(p,q;1-C)$ can be simplified to

\begin{multline*}
\dfrac{1-2C+2C^2}{2}pq-\dfrac12(p+q)+\dfrac{q}{6p}+\dfrac{p}{6q}+(1-2C)\sawtooth{Cpq}+
\dfrac{1}{pq}(\sawtooth{Cpq}^2-\dfrac{1}{12})+\dfrac12-\\
-s(p,q;Cp,0)-s(q,p;Cq,0)-s(p,q;(1-C)p,0)-s(q,p;(1-C)q,0).
\end{multline*}

Hence we prove the following result.

\begin{proposition}\label{prop:siggen}
 If $z=e^{2\pi i C}$ where $C\in[0,1)$ is such that $Cpq$ is not an integer, then the
signature of the torus knot $T_{p,q}$ can be expressed in the following formula.
\begin{multline*}
\sigma(z)=-2(C-C^2)pq+\dfrac{q}{3p}+\dfrac{p}{3q}+(2-4C)\sawtooth{Cpq}+
\dfrac{2}{pq}(\sawtooth{Cpq}^2-\dfrac{1}{12})-\\
-2\left(s(p,q;Cp,0)+s(q,p;Cq,0)+s(p,q;(1-C)p,0)+s(q,p;(1-C)q,0)\right).
\end{multline*}
\end{proposition}

In particular we see rigorously that for large $p$ and $q$ the shape of the function $\sigma(e^{2\pi ix})$
resembles that of the function $2pq(x^2-x)$.

\section{Expressing $\sord(T_{p,q})$ as a rational function}\label{S:proofnogo}

\begin{proposition}\label{th:nogo}
There \emph{does not} exist a rational function $R(p,q)$ such that for all odd and coprime positive integers
\[R(p,q)=\sord(T_{p,q}).\]
\end{proposition}

\begin{proof}
Assume that $R(p,q)=\sigma(T_{p,q})$. Then $S(p,q)=\dfrac14(R(p,q)+(p-1)(q-1))$ is also a rational
function and
\[S(p,q)=\left|\Sigma\cap(0,\dfrac12)\right|=\left|\left\{\dfrac{k}{p}+\dfrac{l}{q}<\dfrac12,\,\,1\le k\le p-1,\,\,1\le l\le q-1.\right\}
\right|\]
(cf. formulae \eqref{eq:S1} and \eqref{eq:S}). If $p|(q-1)$ the value of $S(p,q)$ can be easily computed:

\begin{multline*}
S(p,q)=\sum_{k=1}^{\dfrac{p-1}{2}}\intpart{\dfrac{q}{2}-\dfrac{qk}{p}}=\sum_{k=1}^{\dfrac{p-1}{2}}
\intpart{\dfrac{q-1}{2}-\dfrac{(q-1)k}{p}+\dfrac{p-k}{2p}}=\\
=\sum_{k=1}^{\dfrac{p-1}{2}}\left(\dfrac{q-1}{2}-k\dfrac{q-1}{p}\right)=\dfrac{(q-1)(p-1)^2}{8p}.
\end{multline*}

Since for infinitely many values $(p,q)$ with $q=np+1$ with $p$ odd and $n$ even, we have $p|(q-1)$, it follows
that $S(p,q)=\frac{(q-1)(p-1)^2}{8p}$
on each line $q=np+1$. Since these rational functions agree on infinitely many lines, they must be equal.

But now assume that $p=nq+1$ for some even $n$.
Similar arguments as above show that $S(p,q)$ must also be identical to the function $\frac{(p-1)(q-1)^{2}}{8q}$.
This leads to a contradiction, since these two rational functions are different.
\end{proof}

\begin{remark} We can also compute values of $S(p,q)$ in many other cases, like $q=np-1$, $q=p+2$.
With more care we can prove that e.g.
$S(p,q)-\intfrac{q}{p}$ is not a rational function.
\end{remark}

The proof carries over to show that no such rational function exists for the case $p$ even and $q$ odd. We leave
the obvious details to the reader.

\begin{acknowledgements}
The authors wish to Andrew Ranicki for many remarks during preparation of the paper and  to Julia Collins, Stefan Friedl and
Andr\'as N\'emethi for their interest and comments concerning this article.
\end{acknowledgements}

\end{document}